\documentclass[11pt,leqno]{amsart}

\usepackage[T1]{fontenc}
\usepackage[utf8]{inputenc}
\usepackage{microtype}
\usepackage{mathtools,amssymb,amsthm}
\usepackage{booktabs}
\usepackage{array}
\usepackage{enumitem}
\usepackage[hidelinks]{hyperref}

\numberwithin{equation}{section}
\usepackage{tikz}
\newtheorem{theorem}{Theorem}[section]
\newtheorem{proposition}[theorem]{Proposition}
\newtheorem{lemma}[theorem]{Lemma}
\newtheorem{corollary}[theorem]{Corollary}
\newtheorem{question}[theorem]{Question}

\theoremstyle{definition}

\theoremstyle{remark}
\newtheorem{remark}[theorem]{Remark}
\newcommand{\PP}{\mathbb P}
\newcommand{\CC}{\mathbb C}
\newcommand{\LL}{\mathcal L}
\newcommand{\OO}{\mathcal O}
\newcommand{\Sing}{\operatorname{Sing}}
\newcommand{\mult}{\operatorname{mult}}
\newcommand{\eps}{\varepsilon}

\title[Seshadri constants of small line arrangements]
{Seshadri constants of singular loci of line arrangements with at most twelve lines}
\author[M. Dumnicki]{Marcin Dumnicki}
\address{Faculty of Mathematics and Computer Science, Jagiellonian University, ul. {\L}ojasiewicza 6, 30-348 Krak\'ow, Poland}
\email{Marcin.Dumnicki@uj.edu.pl}

\author[L. Farnik]{{\L}ucja Farnik}
\address{Institute of Mathematics, University of the National Education Commission, Podchor\k{a}\.{z}ych 2, 30-084 Krak\'ow, Poland}
\email{Lucja.Farnik@gmail.com}

\author[K. Hanumanthu]{Krishna Hanumanthu}
\address{Chennai Mathematical Institute, H1 SIPCOT IT Park, Siruseri, Kelambakkam 603103, India}
\email{krishna@cmi.ac.in}

\author[B. Harbourne]{Brian Harbourne}
\address{Department of Mathematics, University of Nebraska--Lincoln, Lincoln, NE 68588, USA}
\email{brianharbourne@unl.edu}

\author[P. Pokora]{Piotr Pokora}
\address{Institute of Mathematics, University of the National Education Commission, Podchor\k{a}\.{z}ych 2, 30-084 Krak\'ow, Poland}
\email{piotr.pokora@uken.krakow.pl}

\author[J. Ro\'e]{Joaquim Ro\'e}
\address{Departament de Matem\`atiques, Universitat Aut\`onoma de Barcelona, 08193 Bellaterra (Barcelona), Spain}
\email{jroe@mat.uab.cat}

\author[H. Tutaj-Gasi\'nska]{Halszka Tutaj-Gasi\'nska}
\address{Faculty of Mathematics and Computer Science, Jagiellonian University, ul. {\L}ojasiewicza 6, 30-348 Krak\'ow, Poland}
\email{halszka.tutaj-gasinska@im.uj.edu.pl}

\date{\today}
\subjclass[2020]{Primary 14C20; Secondary 14N20, 52C30}
\keywords{multipoint Seshadri constants, line arrangements, singular loci, Hesse configuration, incidence structures}

\hypersetup{
  pdftitle={Seshadri constants of singular loci of line arrangements with at most twelve lines},
  pdfauthor={Marcin Dumnicki et al.},
  pdfkeywords={Seshadri constants, line arrangements, singular loci, Hesse configuration}
}

\begin{document}

\begin{abstract}
Let $\LL$ be an arrangement of $d$ lines in the complex projective plane and let
$Z=\Sing(\LL)$.  A question of Pokora asks whether the multipoint Seshadri
constant of $\OO_{\PP^2}(1)$ at $Z$ is always computed by a line of the
arrangement containing the largest number of points of $Z$.  We prove that the
answer is affirmative for every arrangement of at most twelve lines.  The proof
uses a B\'ezout-type covering criterion, elementary restrictions on the
intersection lattice, and a detailed analysis of the two extremal combinatorial
types that can occur for twelve lines.  One type is the Hesse arrangement; it is
handled using a net of quartics through its twelve double points.  The other type
would yield a uniquely determined incidence structure on thirteen lines.  We
show that this incidence structure is not realizable over any field.
\end{abstract}

\maketitle

\section{Introduction}

Seshadri constants measure the local positivity of nef line bundles.  They were
introduced by Demailly \cite{Demailly1992}.  Foundational results and systematic treatments include
\cite{EinLazarsfeld1993,Bauer1999,BauerEtAl2009,Lazarsfeld2004}. For a smooth projective surface $X$, a nef
line bundle $L$, and distinct points $p_1,\ldots,p_r\in X$, the multipoint
Seshadri constant is
\[
 \eps(X,L;p_1,\ldots,p_r)
 =\inf_C\frac{L\cdot C}{\sum_{i=1}^r\mult_{p_i}C},
\]
where the infimum is taken over irreducible reduced curves meeting at least one
of the points.  Even for $X=\PP^2$ and $L=\OO_{\PP^2}(1)$, exact values are
usually difficult to determine.  For very general points, the expected value is
governed by the Nagata conjecture \cite{Nagata1959}.  Configurations arising
from arrangements of curves provide a complementary source of special point
sets for which exact computations can sometimes be made.  The interaction
between line arrangements and algebraic surfaces has a long history; see, for
example, \cite{Hirzebruch1983,OrlikTerao1992}.

Let
\[
 \LL=\{L_1,\ldots,L_d\}\subset \PP^2_{\CC}
\]
be a line arrangement and let
\[
 Z=\Sing(\LL)=\{p\in\PP^2: p\in L_i\cap L_j\text{ for some }i\ne j\}.
\]
We write
\[
 \eps(Z)=\eps\bigl(\PP^2,\OO_{\PP^2}(1);Z\bigr)
\]
and set
\[
 \ell(\LL)=\max_{1\le i\le d}\#(L_i\cap Z).
\]
Since an arrangement line containing $\ell(\LL)$ singular points has Seshadri
ratio $1/\ell(\LL)$, one always has
\[
 \eps(Z)\le \frac{1}{\ell(\LL)}.
\]
The following question was posed in \cite[Question~3.1]{Pokora2019}.

\begin{question}\label{q:main}
Is
\[
 \eps(Z)=\frac{1}{\ell(\LL)}
\]
for every complex line arrangement $\LL$?
\end{question}

Positive answers were obtained in \cite{Pokora2019} for several important
classes, including star arrangements, Fermat arrangements, and the Klein and
Wiman arrangements.  The relation with configurational Seshadri constants was
further studied in \cite{JanaszPokora2020}, and the question was recently placed
in the broader setting of curve configurations on surfaces in
\cite{HanumanthuRoySubramaniam2024}.  The line-arrangement problem remains open
in general.

Our main result settles the first nontrivial range completely.

\begin{theorem}\label{thm:main}
Let $\LL\subset\PP^2_{\CC}$ be an arrangement of $d\le 12$ lines with nonempty
singular locus $Z$.  Then
\[
 \eps(Z)=\frac{1}{\ell(\LL)}.
\]
\end{theorem}

The proof has three ingredients.  First, an effective curve that passes through
$Z$ with a uniform multiplicity gives, by B\'ezout's theorem, a lower bound for
the Seshadri ratios of all curves not contained in it.  Second, the identity
\[
 d-1=\sum_{p\in L_i\cap Z}\bigl(\mult_p\LL-1\bigr)
\]
severely restricts a possible counterexample.  These restrictions settle the
case $d\le 11$.  For $d=12$, only two combinatorial types remain.  One is the
Hesse arrangement.  The other would become, after adjoining a line through six
double points, a thirteen-line incidence structure with six quadruple and
fourteen triple points.  We prove that this structure is unique up to
relabeling and that it has no realization over any field.

Throughout the paper, arrangements are considered over $\CC$, except in
Section~\ref{sec:nonrealizable}, where the nonrealizability statement is valid
over an arbitrary field.

\section{A covering criterion and combinatorial restrictions}

For an irreducible reduced plane curve $D$ meeting $Z$, define its Seshadri ratio
by
\[
 \rho_Z(D)=\frac{\deg D}{\sum_{p\in Z}\mult_pD}.
\]
If $C$ is an effective, possibly nonreduced, plane curve, put
\[
 \mu_Z(C)=\min_{p\in Z}\mult_pC.
\]

\begin{proposition}[Covering criterion]\label{prop:covering}
Let $C$ be an effective plane curve of degree $k$ with $\mu_Z(C)\ge t>0$.
If $D$ is an irreducible reduced curve which is not an irreducible component of
$C$, then
\[
 \rho_Z(D)\ge \frac{t}{k}.
\]
Consequently, if $k/t\le \ell(\LL)$ and every irreducible component $E$ of $C$
satisfies $\rho_Z(E)\ge 1/\ell(\LL)$, then
\[
 \eps(Z)=\frac{1}{\ell(\LL)}.
\]
\end{proposition}

\begin{proof}
Let $e=\deg D$ and $m_p=\mult_pD$.  Since $D$ is not a component of $C$,
B\'ezout's theorem and the local intersection inequality give
\[
 ke=C\cdot D
 \ge \sum_{p\in Z}\mult_pC\,\mult_pD
 \ge t\sum_{p\in Z}m_p.
\]
Thus $\rho_Z(D)\ge t/k$.  Under the additional hypotheses, every irreducible
curve has Seshadri ratio at least $1/\ell(\LL)$, while an arrangement line with
$\ell(\LL)$ singular points gives the opposite inequality.
\end{proof}

The union $L_1+\cdots+L_d$ has multiplicity at least two at every point of $Z$.
We immediately obtain the following useful observation.

\begin{corollary}\label{cor:half}
If $\ell(\LL)\ge d/2$, then
\[
 \eps(Z)=\frac{1}{\ell(\LL)}.
\]
\end{corollary}

\begin{proof}
Apply Proposition~\ref{prop:covering} to $C=L_1+\cdots+L_d$, with $k=d$ and
$t=2$.
\end{proof}

For $p\in Z$, denote by $m_p=\mult_p\LL$ the number of arrangement lines through
$p$, and let
\[
 m(\LL)=\max_{p\in Z}m_p.
\]
For every $L_i\in\LL$, counting the other $d-1$ lines according to their
intersection point with $L_i$ gives
\begin{equation}\label{eq:linecount}
 d-1=\sum_{p\in L_i\cap Z}(m_p-1).
\end{equation}

\begin{lemma}\label{lem:mleell}
If $m(\LL)\ge \ell(\LL)$, then the $m(\LL)$ lines through 
a point of maximal multiplicity cover $Z$, 
and Question~\ref{q:main} has an affirmative answer for $\LL$.

Moreover, if $|Z|>1$ then $m(\LL)\le \ell(\LL)$. 
\end{lemma}

\begin{proof}
If $|Z|=1$ then the claim is obvious, so assume there is no pencil of lines containing all components of $\LL$.
Let $p$ be a point of multiplicity $m=m(\LL)$, and choose an arrangement line
$N$ not passing through $p$.  The line $N$ meets the $m$ lines through $p$ in
$m$ distinct singular points, so $m\le \ell(\LL)$.

Now suppose that $m\ge\ell(\LL)$, and so in fact $m=\ell(\LL)$.  
Every line not through $p$ already contains the
$m$ distinct intersections with the lines through $p$, and hence it contains no
additional singular point.  It follows that every point of $Z$ lies on one of
the $m$ lines through $p$.  Their union has degree $m=\ell(\LL)$ and
multiplicity at least one along $Z$, so Proposition~\ref{prop:covering} applies.
\end{proof}

In particular, the same covering argument applies whenever the arrangement has
a modular point.  Thus it recovers the affirmative answer for supersolvable
line arrangements; see \cite{JambuTerao1984,OrlikTerao1992} for the arrangement-theoretic background.

The preceding results give the basic numerical restrictions used below.

\begin{corollary}\label{cor:restrictions}
If $\LL$ is a counterexample to Question~\ref{q:main}, then
\begin{equation}\label{eq:restrictions}
 m(\LL)<\ell(\LL)<\frac d2,
 \qquad
 \ell(\LL)\bigl(m(\LL)-1\bigr)\ge d-1.
\end{equation}
\end{corollary}

\begin{proof}
The first two inequalities follow from Corollary~\ref{cor:half} and
Lemma~\ref{lem:mleell}.  The last inequality follows from
\eqref{eq:linecount}, since every line contains at most $\ell(\LL)$ singular
points and every summand is at most $m(\LL)-1$.
\end{proof}

\section{Arrangements with at most eleven lines}\label{sec:up-to-eleven}

\begin{proposition}\label{prop:eleven}
Theorem~\ref{thm:main} holds for $d\le 11$.
\end{proposition}

\begin{proof}
Write $\ell=\ell(\LL)$ and $m=m(\LL)$.  A direct check of the integer
solutions of \eqref{eq:restrictions} shows that for $d\le 10$ the only possible
triple is
\[
 (d,\ell,m)=(9,4,3).
\]
In this case equality holds in
$\ell(m-1)\ge d-1$.  Therefore every line contains exactly four triple points
and no other singularities.  The union of the nine arrangement lines has
multiplicity three along $Z$, so Proposition~\ref{prop:covering}, with
$k/t=9/3=3<4$, proves the claim.

Let $d=11$.  The only possible pairs are
\[
 (\ell,m)=(5,3)\quad\text{and}\quad(5,4).
\]
If $m=3$, then \eqref{eq:linecount} forces every line to contain five triple
points.  Counting line--triple-point incidences gives
\[
 55=3t_3,
\]
a contradiction.  Thus this case does not occur.

Assume now that $(\ell,m)=(5,4)$.  If some arrangement line contains no double
point, delete it.  The remaining ten lines still have multiplicity at least two
at every point of $Z$, because every point on the deleted line had multiplicity
at least three.  The resulting degree-ten curve covers $2Z$, and all its
components contain at most five points of $Z$.  Proposition~\ref{prop:covering}
therefore gives $\eps(Z)=1/5$.

We may consequently assume that every line contains a double point.  By
\eqref{eq:linecount}, the list of multiplicities on a line must be one of
\[
 A=(4,4,4,2),\qquad
 B=(4,4,3,2,2),\qquad
 C=(4,3,3,3,2).
\]
Let $a,b,c$ be the numbers of lines of types $A,B,C$, respectively.  Counting
incidences with double and triple points gives
\begin{equation}\label{eq:abc-eleven}
 a+2b+c=2t_2,
 \qquad
 b+3c=3t_3,
 \qquad
 a+b+c=11.
\end{equation}
The second equation implies $3\mid b$, while the first and third imply that $b$
is odd.  Hence $b=3$ or $b=9$.

Suppose first that $b=3$.  Then \eqref{eq:abc-eleven} gives $t_2=7$.  Choose
four auxiliary lines covering the seven double points.  Together with the
arrangement, they form a degree-fifteen curve with multiplicity at least three
at every point of $Z$.  Every arrangement component contains at most five
points of $Z$.  An auxiliary line not belonging to the arrangement also
contains at most five points of $Z$, since otherwise B\'ezout with the union of
the eleven arrangement lines would give
\[
 11\ge 2\cdot 6.
\]
Thus Proposition~\ref{prop:covering} applies with $k/t=15/3=5$.

Finally, suppose that $b=9$.  Then $t_2=10$.  Partition the ten double points
into two sets of five and choose a conic through each set.  Adding these two
conics to the arrangement again gives a degree-fifteen curve with multiplicity
at least three along $Z$.  Its line components have Seshadri ratio at least
$1/5$, as above.  If an irreducible conic component had ratio smaller than
$1/5$, it would contain at least eleven points of $Z$.  Since there are only ten
double points, at least one of these eleven points has multiplicity at least
three in the arrangement.  B\'ezout would then give
\[
 22=2\cdot 11
 \ge \sum_{p\in Q\cap Z}m_p
 \ge 10\cdot 2+3=23,
\]
a contradiction.  Proposition~\ref{prop:covering} completes the proof.
\end{proof}

\section{The case of twelve lines}\label{sec:twelve}

By Corollary~\ref{cor:restrictions}, a counterexample with $d=12$ must satisfy
\[
 \ell(\LL)=5,
 \qquad
 m(\LL)=4.
\]
The identity \eqref{eq:linecount} shows that the possible multiplicity lists on
an arrangement line are
\begin{equation}\label{eq:four-types}
 (4,4,4,3),\quad
 (4,4,4,2,2),\quad
 (4,4,3,3,2),\quad
 (4,3,3,3,3).
\end{equation}

If a line of type $(4,4,4,3)$ occurs, deleting it leaves the same singular set
$Z$, because all four points retain multiplicity at least two.  Some other line
still contains five points, so Proposition~\ref{prop:eleven} gives
$\eps(Z)=1/5$.  If a line of type $(4,3,3,3,3)$ occurs, deletion again leaves
the same singular set.  Let $\ell'$ be the maximum number of points of $Z$ on a
line of the remaining eleven-line arrangement.  Proposition~\ref{prop:eleven}
gives $\eps(Z)=1/\ell'$, whereas the deleted line gives $\eps(Z)\le 1/5$.
Since $\ell'\le 5$, one obtains $\ell'=5$ and again $\eps(Z)=1/5$.

It remains to consider the two line types
\[
 A=(4,4,4,2,2),
 \qquad
 B=(4,4,3,3,2).
\]
Let $a$ and $b$ denote their numbers.  Counting incidences gives
\begin{equation}\label{eq:ab-twelve}
 a+b=12,
 \qquad
 2a+b=2t_2,
 \qquad
 2b=3t_3,
 \qquad
 3a+2b=4t_4.
\end{equation}
The second equation implies that $b$ is even, the third that $3\mid b$, and the
last, after substituting $a=12-b$, that $4\mid b$.  Hence $12\mid b$, so there
are only two cases:
\[
 (a,b)=(12,0)
 \quad\text{or}\quad
 (a,b)=(0,12).
\]

\subsection{The Hesse case}

Assume first that $(a,b)=(12,0)$.  Then
\[
 t_2=12,
 \qquad
 t_3=0,
 \qquad
 t_4=9.
\]
The nine quadruple points and the twelve arrangement lines form the classical
Hesse $(9_4,12_3)$ configuration.  Indeed, each arrangement line contains three
quadruple points, and the twelve lines account for all $\binom{9}{2}=36$ pairs
of quadruple points.  Thus the induced incidence structure is the affine plane
of order three.  Over $\CC$ its line realization is projectively equivalent to
the standard Hesse arrangement; see \cite{ArtebaniDolgachev2009,BassaKisisel2020}.  Let $\omega$ be a primitive cube root of unity.
A convenient equation is
\begin{equation}\label{eq:hesse-arrangement}
 xyz\prod_{i,j=0}^{2}\bigl(x+\omega^i y+\omega^j z\bigr)=0.
\end{equation}
Its twelve double points are
\begin{equation}\label{eq:hesse-double-points}
 (1:0:0),\ (0:1:0),\ (0:0:1),
 \quad\text{and}\quad
 (1:\alpha:\beta),\qquad \alpha^3=\beta^3=1.
\end{equation}
The nine quadruple points are
\begin{equation}\label{eq:hesse-quadruple-points}
 (0:1:-\gamma),\quad (1:0:-\gamma),\quad (1:-\gamma:0),
 \qquad \gamma^3=1.
\end{equation}

Consider the net of quartics
\begin{equation}\label{eq:quartic-net}
 \Lambda=
 \left\langle
 x(y^3-z^3),\;
 y(x^3-z^3),\;
 z(x^3-y^3)
 \right\rangle.
\end{equation}
Its base locus is exactly the set in \eqref{eq:hesse-double-points}.  Indeed, if
$xyz\ne0$, the three equations imply $x^3=y^3=z^3$, while if one coordinate
vanishes, their common zeros are precisely the three coordinate vertices.
Bertini's theorem \cite[Chapter~III, Corollary~10.9]{Hartshorne1977} shows that
a general member is smooth away from the base locus.  At each point
$(1:\alpha:\beta)$ with $\alpha^3=\beta^3=1$, the differential of
$x(y^3-z^3)$ is nonzero, and at each coordinate vertex one of the three
generators has nonzero differential.  Hence singularity at any fixed base point
is a proper linear condition on the coefficients of a member of $\Lambda$; a
general member is therefore smooth at all twelve base points as well.  Finally,
none of the quadruple points is a base point, so a general member of $\Lambda$
avoids all nine of them.

Choose a smooth quartic $Q\in\Lambda$ avoiding the quadruple points.  It is
irreducible and passes simply through the twelve double points.  The divisor
\[
 C=L_1+\cdots+L_{12}+2Q
\]
has degree $20$ and multiplicity four at every point of $Z$.  Every arrangement
line has Seshadri ratio $1/5$, while
\[
 \rho_Z(Q)=\frac{4}{12}=\frac13.
\]
Proposition~\ref{prop:covering} therefore yields
\[
 \eps(Z)=\frac15.
\]

\subsection{The remaining combinatorial type}

Assume now that $(a,b)=(0,12)$.  Then
\[
 t_2=6,
 \qquad
 t_3=8,
 \qquad
 t_4=6.
\]
If the six double points are not collinear, pair them and add the three joining
lines.  The resulting degree-fifteen divisor has multiplicity at least three
along $Z$.  Any line not belonging to the arrangement contains at most six
points of $Z$ by B\'ezout with the twelve arrangement lines.  Equality can occur
only when all six points are double points.  Since the six double points are not
collinear, every component of the covering divisor contains at most five
points of $Z$, and Proposition~\ref{prop:covering} gives $\eps(Z)=1/5$.

It remains to exclude the possibility that the six double points lie on a line
$L_{13}$.  After adjoining $L_{13}$, the six double points become triple
points.  We would obtain a configuration of thirteen lines such that each of
$L_1,\ldots,L_{12}$ contains two quadruple and three triple points, while
$L_{13}$ contains six triple points.  Thus the enlarged configuration would
have
\[
 t_3=14,
 \qquad
 t_4=6.
\]
We now prove that such a configuration does not exist.

\section{A nonrealizable thirteen-line incidence structure}
\label{sec:nonrealizable}

\begin{proposition}\label{prop:unique-incidence}
Up to relabeling, there is a unique incidence structure of thirteen lines with
the following properties:
\begin{enumerate}[label=\textup{(\roman*)}]
 \item each of $L_1,\ldots,L_{12}$ contains two quadruple and three triple
 points;
 \item $L_{13}$ contains six triple points and no quadruple point;
 \item there are six quadruple and fourteen triple points in total.
\end{enumerate}
It is described in Table~\ref{tab:incidence}.
\end{proposition}

\begin{proof}
Use the six quadruple points as vertices of a graph $G$, and represent each of
$L_1,\ldots,L_{12}$ by the edge joining the two quadruple points on that line.
The graph is simple and $4$-regular on six vertices.  Hence its complement is a
perfect matching, and therefore
\[
 G\cong K_6\setminus M_0
\]
for a perfect matching $M_0$.  This determines the quadruple-point incidences uniquely.

\begin{figure}[ht]
\centering
\begin{tikzpicture}[
scale=1.15,
vertex/.style={circle,draw,inner sep=2.5pt},
every node/.style={font=\small}
]
\node[vertex] (P1) at (0,3) {$P_1$};
\node[vertex] (P2) at (3,1.5) {$P_2$};
\node[vertex] (P4) at (3,-1.5) {$P_4$};
\node[vertex] (P6) at (0,-3) {$P_6$};
\node[vertex] (P5) at (-3,-1.5) {$P_5$};
\node[vertex] (P3) at (-3,1.5) {$P_3$};

\draw (P1) -- node[above right] {$L_1$} (P2);
\draw (P1) -- node[above left] {$L_2$} (P3);
\draw (P1) -- node[pos=.42,right] {$L_3$} (P4);
\draw (P1) -- node[pos=.42,left] {$L_4$} (P5);

\draw (P2) -- node[pos=.58,right] {$L_5$} (P6);
\draw (P3) -- node[pos=.58,left] {$L_6$} (P6);
\draw (P4) -- node[below right] {$L_7$} (P6);
\draw (P5) -- node[below left] {$L_8$} (P6);

\draw (P2) -- node[above] {$L_9$} (P3);
\draw (P2) -- node[right] {$L_{10}$} (P4);
\draw (P3) -- node[left] {$L_{11}$} (P5);
\draw (P4) -- node[below] {$L_{12}$} (P5);

\draw[densely dashed] (P1) -- node[pos=.42,right] {$P_1P_6$} (P6);
\draw[densely dashed] (P2) -- node[pos=.42,below] {$P_2P_5$} (P5);
\draw[densely dashed] (P3) -- node[pos=.42,above] {$P_3P_4$} (P4);

\end{tikzpicture}
\caption{The graph $G=K_6\setminus M_0$,
with the edges labeled according to
Table~\ref{tab:incidence}. The dashed edges form the complementary
perfect matching $M_0={P_1P_6,P_2P_5,P_3P_4}$.}
\label{fig}
\end{figure}
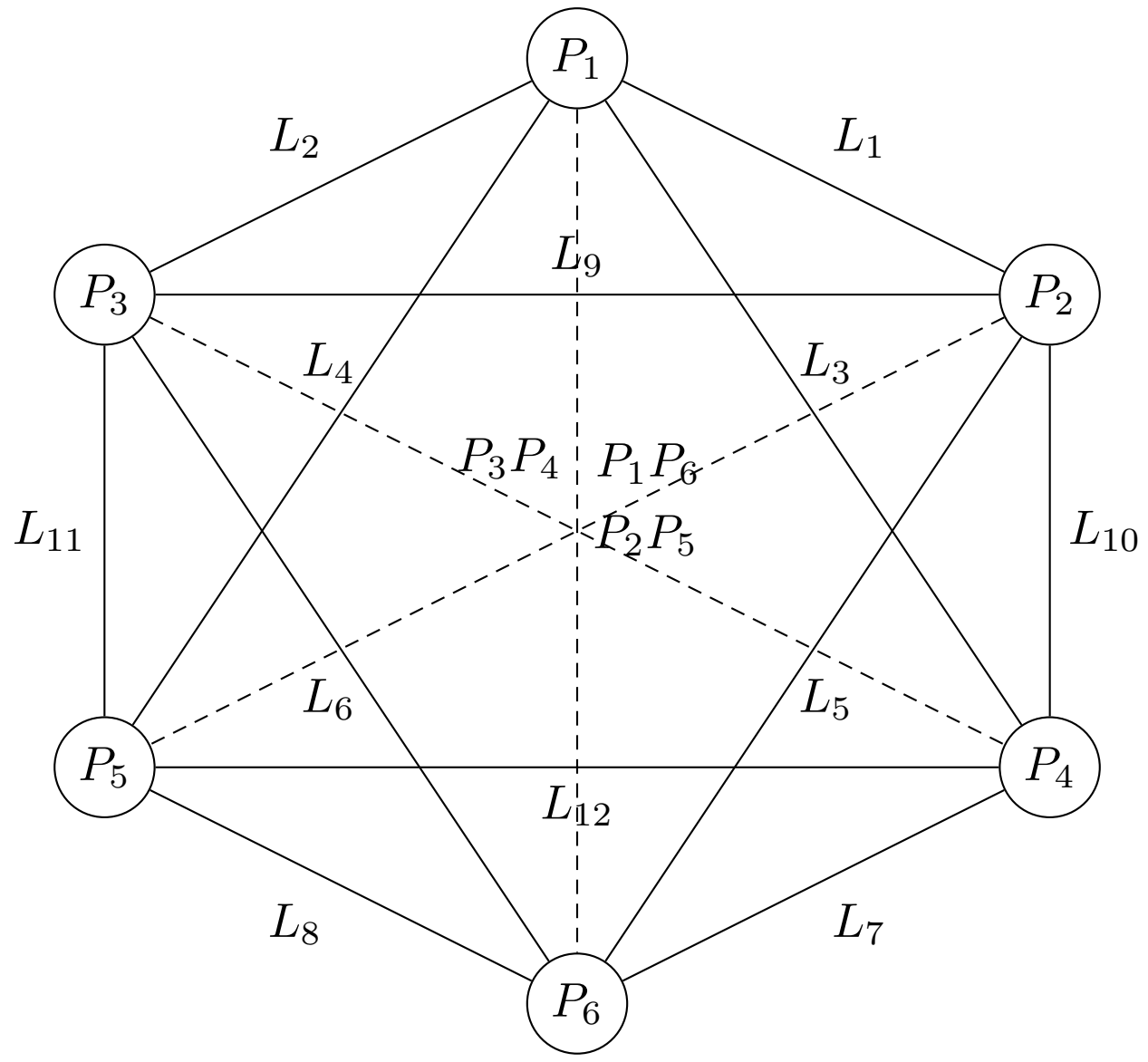

Two of the first twelve lines which do not share a quadruple point correspond
to disjoint edges of $G$.  A triple point not on $L_{13}$ corresponds to a
perfect matching of $G$.  The graph $G$ has exactly eight perfect matchings. Indeed, fix the vertex $P_1$. It can be matched with any of its four neighbours in $G$. After
such a choice, the subgraph induced by the four remaining vertices is
$K_4$ with exactly one edge removed, and hence has precisely two perfect
matchings. Consequently, the total number of perfect matchings of $G$ is
\[
4\cdot 2=8.
\]
There are exactly eight triple points away from $L_{13}$, so all these
perfect matchings occur. A pair of disjoint edges of $G$ is contained in
such a perfect matching unless the two unused vertices form an edge of
$M_0$. The six remaining pairs are disjoint and exhaust the twelve edges
of $G$; they are therefore precisely the six pairs meeting on
$L_{13}$. This determines all incidences and gives
Table~\ref{tab:incidence}.
\end{proof}

\begin{table}[ht]
\centering
\caption{The unique incidence structure from Proposition~\ref{prop:unique-incidence}.}
\label{tab:incidence}
\small
\begin{tabular}{c@{\qquad}l|c@{\qquad}l}
\toprule
point & incident lines & point & incident lines\\
\midrule
$P_1$  & $L_1,L_2,L_3,L_4$       & $P_{11}$ & $L_2,L_7,L_{13}$\\
$P_2$  & $L_1,L_5,L_9,L_{10}$    & $P_{12}$ & $L_2,L_8,L_{10}$\\
$P_3$  & $L_2,L_6,L_9,L_{11}$    & $P_{13}$ & $L_3,L_5,L_{11}$\\
$P_4$  & $L_3,L_7,L_{10},L_{12}$ & $P_{14}$ & $L_3,L_6,L_{13}$\\
$P_5$  & $L_4,L_8,L_{11},L_{12}$ & $P_{15}$ & $L_3,L_8,L_9$\\
$P_6$  & $L_5,L_6,L_7,L_8$       & $P_{16}$ & $L_4,L_5,L_{13}$\\
$P_7$  & $L_1,L_6,L_{12}$        & $P_{17}$ & $L_4,L_6,L_{10}$\\
$P_8$  & $L_1,L_7,L_{11}$        & $P_{18}$ & $L_4,L_7,L_9$\\
$P_9$  & $L_1,L_8,L_{13}$        & $P_{19}$ & $L_9,L_{12},L_{13}$\\
$P_{10}$ & $L_2,L_5,L_{12}$      & $P_{20}$ & $L_{10},L_{11},L_{13}$\\
\bottomrule
\end{tabular}
\end{table}

\begin{theorem}\label{thm:nonrealizable}
The incidence structure in Table~\ref{tab:incidence} is not realizable by lines in $\PP^2$ over any field.
\end{theorem}
\begin{remark}
The incidence table defines a genuine abstract incidence structure, which is unique up to relabeling. Thus the obstruction is not combinatorial consistency but representability. Over $\mathbb C$, its nonrealizability can already be detected by the Hirzebruch-type inequality \cite[Remark 11]{PokH}
\[
 2t_2+t_3\geq d+3+\sum_{r\geq 5}(r-4)t_r.
\]
Indeed, here $d=13$, $t_2=0$, and $t_3=14$, so that $14<16$.
The coordinate argument below is stronger, since it excludes realizations over
an arbitrary field.
\end{remark}
\begin{proof}
Assume that a realization exists over a field $K$.  The four points
$P_1,P_9,P_{16},P_6$ are in general position, so after a projective change of
coordinates we may take
\[
 \begin{aligned}
 P_1&=(1:0:0), & P_9&=(0:0:1),\\
 P_{16}&=(0:1:0), & P_6&=(1:1:1).
 \end{aligned}
\]
The incidence table then gives
\[
 L_{13}:x=0,\quad L_1:y=0,\quad L_4:z=0,
 \quad L_5:x-z=0,\quad L_8:x-y=0.
\]
Write
\[
 P_{14}=(0:t:1),\qquad P_{11}=(0:s:1).
\]
The points are distinct, so $s\ne t$.  The lines through the prescribed pairs
are
\[
 \begin{aligned}
 L_3&:y-tz=0,&
 L_6&:(1-t)x-y+tz=0,\\
 L_2&:y-sz=0,&
 L_7&:(1-s)x-y+sz=0.
 \end{aligned}
\]
A direct intersection computation yields
\[
 \begin{aligned}
 P_2&=(1:0:1),&
 P_3&=(t-s:s(t-1):t-1),\\
 P_4&=(t-s:t(1-s):1-s),&
 P_{15}&=(t:t:1),\\
 P_{18}&=(-1:s-1:0),&
 P_{12}&=(s:s:1),\\
 P_{17}&=(-1:t-1:0).&&
 \end{aligned}
\]
Consequently, the lines $L_9=P_2P_3$ and $L_{10}=P_2P_4$ have equations
\[
 \begin{aligned}
 L_9&:s(t-1)x+(s-1)y-s(t-1)z=0,\\
 L_{10}&:t(s-1)x+(t-1)y-t(s-1)z=0.
 \end{aligned}
\]
The required incidences $P_{15},P_{18}\in L_9$ and
$P_{12},P_{17}\in L_{10}$ give
\begin{align}
 f_1&=st^2-st+s-t=0,\label{eq:f1}\\
 f_2&=s^2-st-s+1=0,\label{eq:f2}\\
 f_3&=s^2t-st-s+t=0,\label{eq:f3}\\
 f_4&=t^2-st-t+1=0.\label{eq:f4}
\end{align}
Subtracting \eqref{eq:f3} from \eqref{eq:f1} and \eqref{eq:f4} from
\eqref{eq:f2} gives
\[
 (s-t)(st-2)=0,
 \qquad
 (s-t)(s+t-1)=0.  
\]
Since $s\ne t$, we obtain $st=2$ and $s+t=1$.  Substitution into
\eqref{eq:f1} gives
\[
 f_1=2t-2+s-t=s+t-2=-1,
\]
a contradiction in every field.
\end{proof}

\begin{proof}[Proof of Theorem~\ref{thm:main}]
Proposition~\ref{prop:eleven} treats $d\le 11$.  Section~\ref{sec:twelve}
reduces the case $d=12$ to the Hesse case and the possible collinearity of six
double points.  The Hesse case was settled using the quartic covering divisor,
and Theorem~\ref{thm:nonrealizable} excludes the remaining incidence structure.
\end{proof}

\section{A counterexample with sixteen lines over $\overline{\mathbb{F}}_7$}
\label{sec:char7}

Question~\ref{q:main} is stated over $\CC$, but the definition of $\eps(Z)$
makes sense over any algebraically closed field.  The results of Section~2
remain valid in that generality, since their proofs use only B\'ezout's theorem
and the identity~\eqref{eq:linecount}.  It is therefore natural to ask whether
Question~\ref{q:main} has an affirmative answer in positive characteristic.  In
this section we show that it does not.  Over $\overline{\mathbb{F}}_7$ there is
an arrangement of sixteen lines, each containing at most seven singular points,
together with a line outside the arrangement that contains eight.  We
then show that this phenomenon is confined to characteristic~$7$, i.e., the
intersection lattice of the arrangement has no realization in any other
characteristic.

\subsection{The arrangement}
\label{subsec:arr}
Let $\LL=\{L_1,\ldots,L_{16}\}$ be the arrangement in
$\PP^2_{\overline{\mathbb{F}}_7}$ whose lines are the zero loci of the
following linear forms:
\begin{equation}\label{eq:sixteen-lines}
\arraycolsep=4pt
\begin{array}{llll}
 L_1\colon y+2z, & L_2\colon y+3z, & L_3\colon y+5z, & L_4\colon x,\\[2pt]
 L_5\colon x+z, & L_6\colon x+3z, & L_7\colon x+y+z, & L_8\colon x+y+6z,\\[2pt]
 L_9\colon x+2y+3z, & L_{10}\colon x+3y+4z, & L_{11}\colon x+4y+z, & L_{12}\colon x+4y+5z,\\[2pt]
 L_{13}\colon x+5y+z, & L_{14}\colon x+5y+6z, & L_{15}\colon x+6y+z, & L_{16}\colon x+6y+2z.
\end{array}
\end{equation}
All sixteen lines are defined over $\mathbb{F}_7$, and hence so is every point of
$Z=\Sing(\LL)$.  A direct computation shows that $|Z|=38$, with
\[
 t_2=16,\qquad t_3=12,\qquad t_4=8,\qquad t_5=2.
\]
The two quintuple points are $(1:0:6)$ and $(1:1:6)$.  Fourteen of the lines contain
seven points of $Z$, and the remaining two, $L_{10}$ and $L_{15}$, contain six.
Hence
\[
 \ell(\LL)=7,\qquad m(\LL)=5.
\]

Now consider the line
\[
 M\colon x+6y+6z=0,
\]
which does not belong to $\LL$.  It contains the following eight points of $Z$,
each of which is a double point of $\LL$:
\begin{center}
\small
\begin{tabular}{cc@{\qquad}cc}
\toprule
point & lines through it & point & lines through it\\
\midrule
$(1:2:6)$ & $L_1,L_5$    & $(1:6:2)$ & $L_6,L_{12}$\\
$(1:5:3)$ & $L_2,L_{10}$ & $(1:0:1)$ & $L_8,L_{14}$\\
$(1:3:5)$ & $L_3,L_{13}$ & $(1:4:4)$ & $L_9,L_{11}$\\
$(0:1:6)$ & $L_4,L_7$    & $(1:1:0)$ & $L_{15},L_{16}$\\
\bottomrule
\end{tabular}
\end{center}
Thus $\rho_Z(M)=1/8$, and therefore
\[
 \eps(Z)\le\frac18<\frac17=\frac1{\ell(\LL)}.
\]
So $\LL$ is a counterexample to Question~\ref{q:main} over
$\overline{\mathbb{F}}_7$.

\begin{remark}\label{rem:extremal-line}
The line $M$ is extremal.  A line $N\notin\LL$ meets each arrangement line in
exactly one point, so
\[
 2\,\#(N\cap Z)\le\sum_{p\in N\cap Z}m_p\le d.
\]
Hence $N$ contains at most $d/2$ points of $Z$, with equality only if $N$
passes through $d/2$ double points and the arrangement lines meet $N$ in pairs.
The table above shows that $M$ attains this bound for $d=16$.  Exactly the same
situation arose for $d=12$ in Section~\ref{sec:nonrealizable}, where six double
points on a line $L_{13}$ would have given $\eps(Z)\le 1/6<1/5$.
Theorem~\ref{thm:nonrealizable} rules that out over every field, whereas for
$d=16$ it occurs in characteristic~$7$.  As it must, the triple
$(d,\ell,m)=(16,7,5)$ satisfies the restrictions~\eqref{eq:restrictions}.
\end{remark}

\begin{remark}
This arrangement was found by one of the authors with the help of an AI system
(the large language model Claude Opus~5).  Its correctness was then checked
carefully by the authors.
\end{remark}

\subsection{No lift to characteristic zero}

The counterexample is a genuinely positive-characteristic phenomenon, as the
following proposition shows.  We keep the labeling of~\eqref{eq:sixteen-lines}.

\begin{proposition}\label{prop:char7}
Let $K$ be a field, and let $\LL'$ be an arrangement of sixteen lines in
$\PP^2_K$ with the same intersection lattice as $\LL$.  Then $K$ has
characteristic~$7$.  In particular, $\LL$ has no complex realization.
\end{proposition}

\begin{proof}
No three of $L_4,L_8,L_9,L_{12}$ are concurrent.  After a projective change of
coordinates we may therefore assume that
\[
 L_4\colon x=0,\qquad L_8\colon y=0,\qquad L_{12}\colon z=0,\qquad
 L_9\colon x+y+z=0.
\]
The concurrences of $\LL$ now determine seven further lines, one at a time.
Each line in the following table passes through two intersection points of
lines that are already known.  In every row these two points are distinct over
any field, so they determine the line.
\begin{center}
\small
\begin{tabular}{lll}
\toprule
line & passes through & equation\\
\midrule
$L_{10}$ & $L_9\cap L_{12}=(1:-1:0)$,\ \ $L_4\cap L_8=(0:0:1)$ & $x+y=0$\\
$L_{13}$ & $L_4\cap L_{12}=(0:1:0)$,\ \ $L_8\cap L_9=(1:0:-1)$ & $x+z=0$\\
$L_{16}$ & $L_8\cap L_{12}=(1:0:0)$,\ \ $L_4\cap L_9=(0:1:-1)$ & $y+z=0$\\
$L_1$    & $L_8\cap L_{12}=(1:0:0)$,\ \ $L_{10}\cap L_{13}=(1:-1:-1)$ & $y-z=0$\\
$L_{11}$ & $L_1\cap L_4=(0:1:1)$,\ \ $L_{10}\cap L_{16}=(1:-1:1)$ & $2x+y-z=0$\\
$L_5$    & $L_9\cap L_{12}=(1:-1:0)$,\ \ $L_{11}\cap L_{13}=(1:-3:-1)$ & $x+y-2z=0$\\
$L_6$    & $L_{13}\cap L_{16}=(1:1:-1)$,\ \ $L_8\cap L_{11}=(1:0:2)$ & $2x-3y-z=0$\\
\bottomrule
\end{tabular}
\end{center}
Finally, $L_4$, $L_5$ and $L_6$ are concurrent in $\LL$; in the coordinates
of~\eqref{eq:sixteen-lines} all three pass through $(0:1:0)$.  In the present
coordinates, $L_4\cap L_5=(0:2:1)$, and substituting this point into the
equation of $L_6$ gives
\[
 2\cdot 0-3\cdot 2-1=-7.
\]
Hence $-7=0$ in $K$, that is, $\operatorname{char}K=7$.
\end{proof}

\section{Further questions}

Theorem~\ref{thm:main} shows that the smallest counterexample to Question~\ref{q:main} requires at least thirteen lines.  Corollary~\ref{cor:restrictions}
continues to impose
\[
 m(\LL)<\ell(\LL)<\frac d2,
 \qquad
 \ell(\LL)(m(\LL)-1)\ge d-1.
\]
The covering criterion suggests that the most promising candidates are
arrangements for which the singular points of smallest multiplicity cannot be
covered efficiently by low-degree curves.  It would be interesting to combine
these restrictions with databases of realizable rank-three matroids and with
interpolation computations for fat point schemes.

\begin{question}
Does Question~\ref{q:main} have an affirmative answer for arrangements of
thirteen lines?  More generally, is there any complex line arrangement for
which the Seshadri constant of the full singular locus is not computed by an
arrangement line?
\end{question}
\section*{Acknowledgments}
KH is supported by 
a grant from Infosys Foundation and 
ANRF MATRICS grant (ANRF/ARGM/2025/002374/MTR). 
The project originated during the mini-workshop \emph{Seshadri Constants}, held at the Mathematisches Forschungsinstitut Oberwolfach from 10 to 16 November 2019.
\section*{AI Declaration}
Apart from the example in Section \ref{subsec:arr}, which was found by \texttt{Claude}, no AI tools were used in preparing this manuscript.
\bibliographystyle{amsplain}
\bibliography{references}

\end{document}